\magnification=\magstep1
\documentstyle{amsppt}
\NoRunningHeads
\font\tyt=cmbx12

\def\Om{\Omega }
\def\MA{Monge-Amp\`ere }
\def\psh{plurisubharmonic }
\def\cl{\centerline }
\def\nb{neighbourhood }
\def\fii{\varphi }
\def\we{\wedge }

\def\om{\omega }
\def\db{\overline{\partial } }
\def\dd{\partial }
\def\ep{\epsilon }

\def\RHS{right hand side}
\def\de{\delta }
\def\ze{\zeta }
\def\al{\alpha }

\def\du{(dd^c u )^n }

\def\cn{\Bbb C ^n }

\def\MAE{Monge-Amp\`ere equation}

\def\po{\partial\Omega }

\def\ude{u_{\de } }

\bigskip
\def\sk{\nopagebreak\vskip 3mm}

\phantom{a}\vskip 25mm

\cl{\tyt H\"older continuity of solutions to  the complex}

\sk

\cl{\tyt  \MAE \ with the \RHS\  in $L^p .$}

\sk

\cl{\tyt The case of compact K\"ahler manifolds}

\sk
\document

\topmatter
\author S\l awomir Ko\l odziej \endauthor
\affil Jagiellonian University, Cracow \endaffil

\bigskip

\abstract We prove that on compact K\"ahler manifolds solutions to
the complex  \MAE , with the \RHS\  in $L^p ,\  p>1,$  are
H\"older continuous.
\endabstract

\subjclass
primary 32U05, secondary 32U40
\endsubjclass
\keywords
plurisubharmonic function, complex Monge-Amp\`ere operator
\endkeywords
\endtopmatter

 Let $M$ be a compact n-dimensional K\"ahler manifold with the
 fundamental form $\om $ given in local coordinates by
$$
\om =\frac{i}{2}\sum_{k,j}g_{k\bar j } dz^k \wedge d\bar z ^j .
 $$

An upper semicontinuous function $u$ on $M$ is called
$\om$-plurisubharmonic ($\om$-psh \rm in short) if $ dd^c u +\om
\geq 0.$

We study the regularity of  $\om$-psh solutions $u $ of the
complex \MAE
$$
(\om +dd^c u  )^n = f\om ^n , \tag{0.1} $$ where $f\in L^1 (M),
f\geq 0, \int _M f\om ^n = \int _M \om ^n  $ is a given function.
For smooth positive $f$ the equation was solved by Yau [Y]. Later
it was shown in [K1] that for $f\in L^p (M), \ p>1$ there exists a
continuous solution. From [K2] we know that $L^{\infty }$ norm of
a difference of (suitably normalized) solutions is controlled by
$L^1$ norm of the difference of functions on the \RHS\ (see
Theorem 1.1 below). Here we shall prove that for $f\in L^p (M), \
p>1$ the solutions are H\"older continuous. The exponent depends
on $M$ and $||f||_p$.

The corresponding result in strictly pseudoconvex domains has
been obtained in [GKZ].

The results are motivated, in part, by a recent work of Tian and
Zhang [TZ] where the authors study the K\"ahler-Ricci flow on
projective manifolds. Later more papers on the subject appeared in
ArXiv: [ST],[EGZ],[Z].  When the canonical divisor is big and nef
the flow initiated at any K\"ahler metric tends to a current which
is a smooth K\"ahler-Einstein metric off a subvariety $S$ of $M$.
The potential of this current is continuous also along the
singular set. Since the \RHS\ of the \MAE\ satisfied by the
potential blows up along $S$ at the rate $d^{\al }$, where $d$
denotes the distance from $S$, and it is integrable at the same
time, it belongs to some $L^p , p>1 .$ The K\"ahler form on the
left hand side also degenerates, so   the \MAE\  here becomes more
complicated.  However,  our result, with the same proof, holds on
compact K\"ahler orbifolds. Also, if the limit metric has
singularities that can be blown down, then we may pull-back the
equation from a K\"ahler manifold. Thus Theorem 2.1
 can be applied to some 2 - dimensional examples of
the K\"ahler-Ricci flow considered in [TZ] and [ST].

I would like to thank G. Tian for an invitation to Princeton and
the possibility of discussing this topic. I also thank V. Guedj,
M. Paun and Z. Zhang for their comments on this paper.  Z. Zhang's
suggestion  shortened the proof.

\bigskip
\bf 1. Preliminaries. \rm
\bigskip

For the background of the definitions and results that follow see
[K3]. Using the differential operators
 $d=\dd +\db ,\ d^c =i(\db -\dd
)$ we define for given  bounded \psh function $u$  the \MA\
measure
$$
\du = dd^c u \we dd^c u \we ... \we dd^c u \ \ \ (n \text{ terms})
$$
(see [BT1]). This is a nonnegative Borel measure.

Let us recall a stability statement from [K2] concerning the
equation \thetag{0.1}. On a compact K\"ahler manifold $M$ with a
fundamental form $\om$ the $L^p$ norms are defined by
$$
||f||_p =(\int _M |f|^p \om ^n )^{1/p} .
$$

\proclaim{Theorem 1.1} Given $p>1$, $\ep >0, c_0 >0$ and $||f||_{p
} <c_0 , ||g||_{p } <c_0$ satisfying the normalizing condition in
\thetag{0.1} there exists $c(\ep , c_0 )$ with
$$
||\fii -\psi ||_{\infty } \leq c(\ep , c_0 )||f-g||_1 ^{1/(n+3+\ep
)}
$$
for suitably normalized $\fii $ and $\psi $.
\endproclaim
Let $\Om $ be a domain in $\cn $ and $u\in PSH(\Om ).$
 For $z\in \Om _{\de }: =\{ z\in \Om : dist (z, \po )>\de \}$
define  a \psh function

$$\tilde{u}_{\de }(z) =[\tau (n)\de ^{2n}]^{-1} \int _{|\ze |\leq \de  }
u(z+\ze )\, dV(\ze ), \ \ \tau (n):=\int _{|\ze |\leq 1 } \,
dV(\ze ) ,
$$
where $dV$ denotes the Lebesgue measure.  In [GKZ] it is proved
that
$$
\int _{ \Om _{\de } } (\tilde{u}_{\de } - u)\, dV(\ze )\leq c_1
(n)||\Delta u||_1 \de ^2 ,\tag{1.1}
$$
with the constant $c_0 (n)$ depending only on the dimension.  For
the sake of completeness we include the proof here. Applying
Jensen's formula and Fubini's theorem we obtain the following
estimates (with $ \sigma _{2n-1}$ denoting the surface measure of
the unit sphere)
$$\aligned
&\int _{ \Om _{\de } } \de ^{-2}(\tilde{u}_{\de } - u) (\ze )\,
dV(\ze )
\\
=&\frac{2n}{\de ^{2(n+1) } \sigma _{ 2n-1} }\int _{ \Om _{\de }
}\int _0 ^{\de } r^{2n-1} \int _0 ^r t^{1-2n} \int _ {|\ze
|\leq t } \Delta u(z+\ze )\, dV(\ze ) \, dt\, dr\, dV(z) \\
=& \frac{2n}{\de ^{2(n+1) } \sigma _{ 2n-1} }\int _0 ^{\de }
r^{2n-1} \int _0 ^r t^{1-2n}  \int _ {|\ze
|\leq t } \int _{ \Om _{\de } }\Delta u(z+\ze )\, dV(z)\, dV(\ze ) \, dt\, dr \\
\leq & \frac{2n}{\de ^{2(n+1) } \sigma _{ 2n-1} }\int _0 ^{\de }
r^{2n-1} \int _0 ^r t^{1-2n} \int _ {|\ze |\leq t } ||\Delta
u||_1\, dV(\ze ) \, dt\, dr \\ = & c_0 (n)||\Delta u||_1 ,
\endaligned
$$
where $||\Delta u||_1 = \int _ {\Om } \Delta u(\ze )\, dV(\ze ) $.
The estimate \thetag{1.1} directly follows from this one.

\bigskip
\bf 2. Main theorem. \rm
\bigskip

\proclaim{Theorem 2.1} For $p>1$ and $f \in L^p (M) $ satisfying
the normalizing condition in \thetag{0.1} the solution $u$ of
\thetag{0.1} is  H\"older continuous with the H\"older exponent
which depends on $ p$,  $M$ and $||f||_p$.
\endproclaim
\demo{Proof} It follows from [K1] that $||u||_{\infty }$ depends
on $ p$,  $M$ and $||f||_p$. We can assume that $1<u .$ Choose a
finite number of coordinate balls $B''_j =B(a_j ,3r )$ such that
$B'_j = B(a_j , r)$ cover $M$ and denote by $B_j$ the balls $B(a_j
,2r )$. Since the transition functions for those charts have
bounded Jacobians one can find a constant $C>0$ depending only on
$M$ and such that for all $\de <r/2C$ and $z\in B_j \cap B_k$ we
have
$$
B_j (z ,\de )\subset B_k (z, \frac{C}{2}\de ) ,\tag{2.1}
$$
where $B_j (z ,\de )$ denotes the ball centered at $z$ of radius
$\de$ in the chart $B''_j.$

Fix non positive functions $\rho _j \in C^{\infty } (B_j )$ with
$\rho _j =0$ on $B'_j$, $-1\leq \rho _j \leq 0$ in $B_j$, and
$\rho _j =-1$ on a \nb of $\partial B_j .$ For some $c_1$
depending on $M$ we have
$$
dd^c \rho _j \geq -c_1 \om .
$$
For fixed $\ep >0$ we choose $N>2C$  and $\al < \frac{1}{q(n+3+\ep
)+1 }$ (with $p,q$ conjugate) satisfying
$$
2(2c_1 ||u||_{\infty } +1)<N^{-\al }\frac{\log N}{\log C}
.\tag{2.2}
$$
It is possible since we can choose $N$ so big that $\frac{\log
N}{\log C}$ is bigger than two times the left hand side, and then
we take $\al$ small enough. Using the coordinates in $B''_j$ we
define regularizations
$$u_{j,\de }(z) =\max _{|w|<\de } u(z+w), \ \ \ z\in B_j .
$$
Let us yet define two auxiliary functions
$$
\chi (\de )=\de ^{-\al } \max _j \max _{z\in B_j } (u_{j,\de }
-u)(z) ,
$$
and
$$
\eta (\de )= \max _j \max _{z\in B_j } (u_{j,C\de } -u_{j,\de
})(z) ,
$$
By \thetag{2.1}   we have
$$
\max _{z\in B_j \cap B_k} |(u_{j,\de } - u_{k,\de } )(z)|\leq \eta
(\de ) .\tag{2.3}
$$
We shall approximate $u$ by $\om$-psh functions $\ude$ which are
created by gluing together the local regularizations $u_{j,\de }$
(comp. [D]). The function $\eta$ defined above measures the
correction term when we pass from local to global regularization.
Note that, by continuity of $u$, $$\lim _{\de \to 0} \eta (\de )
=0.$$ Set
$$
\ude (z)= (1+C_1 \eta (\de ) )^{-1} \max _j (u_{j,\de }(z) + \eta
(\de )\rho _j (z) ), \ \ \ C_1 =2c_1 .
$$
It is continuous on $M$ since, by \thetag{2.3}, the maximum on the
\RHS\  must be attained for $j$ such that $z\in B'_j .$  Note also
that since for $c_1 \eta (\de ) <1$
$$ dd^c (u_{j,\de }(z) + \eta (\de )\rho _j (z) )\geq
-(1+\frac{C_1}{2}\eta (\de ))\om$$ one obtains, via an inequality
from [BT1] estimating $dd^c \max (u,v)$ from below,
$$
dd^c \ude +\om >0 ,\tag{2.4}
$$
if $\de$ is sufficiently small. To finish the proof we need to
verify the following claim.
\bigskip

\bf Claim. \rm $\chi $ is bounded on some nonempty interval $(0,
\tilde{\de }  ).$
\bigskip

Suppose that $\chi(\de )
>\max (9, \chi (N\de ))$ and $N\de <r/2 .$ Then the set
$$
E=\cup _j \{ z\in B_j :  (u_{j,\de } -u)(z) > (\frac{\chi (\de
)}{3} -2)\de ^{\al } \}
$$
is nonempty. Take $g=0$ on $E$ and $g=C_2 f$ on $M\setminus E$
with the constant $C_2$ chosen so that $\int _M g\om ^n =\int _M
\om ^n .$

Now we compare $u_{j,\de }$ with
$$\tilde{u}_{j, \de }(z) =[\tau (n)\de ^{2n}]^{-1} \int _{|\ze |\leq \de  }
u(z+\ze )\, dV(\ze ), \ \ \tau (n):=\int _{|\ze |\leq 1 } \,
dV(\ze ) ,
$$
where the coordinates of $B''_j$ are used. Given $z\in B_j$ we
find $w_z$ with $|w_z|=\de $ such that
$$\aligned
u_{j,\de } &= u(z+w_z )  \leq \tilde{u}_{j, \sqrt{\de }}(z+w_z )
\leq \tilde{u}_{j, \sqrt{\de }}(z) +2||u||_{\infty } \sqrt{\de } .
\endaligned
$$
(Note that defining $\tilde{u}_{j, \sqrt{\de }}(z)$ and
$\tilde{u}_{j, \sqrt{\de }}(z+w_z )$ we integrate over the same
domain except for the piece of volume at most  $2\tau (n)\de ^{
n+\frac{1}{2} }.$)

Since $\al <1/2$ we infer from this estimate for $\de <\de _0$ and
$\de _0$ small enough
$$
E\cap B_j \subset \{ u_{j,\de } -u > \de ^{\al } \} \subset \{
\tilde{u}_{j, \sqrt{\de }}-u
> \de ^{\al }/2 \} .
$$
Thus, as $||\Delta u||_1$ is a priori bounded on every $B''_j$,
applying \thetag{1.1} one obtains
$$
\int_{E\cap B_j } \om ^n < c_3 \de ^{1- \al }$$ for all $j$ and
consequently
$$
\int _E \om ^n <c_4 \de ^{1- \al } ,$$ with the constant depending
only on $M$. Hence, upon the use of H\"older inequality
$$
\int _E f\om ^n \leq ||f||_p (\int _E \om ^n )^{1/q} \leq c_5 \de
^{(1- \al )/q } ,$$ where $c_5$ depends also on $||f||_p $. By
Theorem 1.1 for $\de <\de _1$ and $\de _1$ small enough if $v$
solves
$$(\om +dd^c v )^n =g\om ^n$$ and is suitably normalized then
$$
||u-v||_{\infty } \leq ||f-g||^{\frac{1}{(n+3+\ep )} } \leq c_6
\de ^{\frac{1-\al }{q(n+3+\ep ) } } \leq \de ^{\al } \tag{2.5}
$$
(since by our choice of $\al$ we have $\al < \frac{1-\al
}{q(n+3+\ep ) } $).

\proclaim{Proposition} If we choose $z_0\in B_{j_0}$ so that
$$
 (u_{j_0 ,\de } -u)(z_0 ) =\chi (\de )\de ^{\al } ,
$$
then
$$
\sup _{M\setminus E}(\ude -v)<  (\ude -v)(z_0 ).
$$
\endproclaim
\demo{Proof} Take $z\in (M\setminus E )\cap B_j .$ Then
$$
 (u_{j,\de } -u)(z )\leq (\frac{\chi (\de
)}{3} -2)\de ^{\al }
$$
and therefore, by \thetag{2.5}
$$
 (u_{j,\de } -v)(z )\leq (\frac{\chi (\de
)}{3} -1)\de ^{\al }  .
$$
Since $u>1$ we get from this
$$
(\ude -v)(z)\leq \max _{j: z\in B_j } (u_{j,\de } -v)(z )\leq
(\frac{\chi (\de )}{3} -1)\de ^{\al }  .\tag{2.6}
$$
Again,  by \thetag{2.5}
$$
(u_{j_0 ,\de } -v)(z_0 ) \geq (\chi (\de )-1)\de ^{\al } .
$$
Therefore the definition of $\ude$ yields
$$
(\ude -v)(z_0 )\geq (\chi (\de )-1)\de ^{\al } -\eta (\de ) (2c_1
||u||_{\infty } +1 ) .\tag{2.7}
$$
The Three Circles Theorem gives for $\de$ small enough
$$
(u_{j,N\de } -u_{j,\de })\geq \frac{\log N}{\log C}(u_{j,C\de }
-u_{j,\de }).
$$
It follows that, choosing $j$ so that
$$
\eta (\de )=  \max _{z\in B_j } (u_{j,C\de } -u_{j,\de })(z)
$$
we obtain
$$
(N\de )^{\al } \chi (N\de )\geq \frac{\log N}{\log C}\eta (\de ).
$$
Further, since $\chi (\de )\geq \chi (N\de )$, we get from
\thetag{2.2} that
$$
\de ^{\al }\chi (\de ) \geq \frac{\log N}{\log C}\eta (\de
)N^{-\al } > 2\eta (\de ) (2c_1 ||u||_{\infty }   +1 ). $$
Inserting this into \thetag{2.7} we finally arrive at
$$
(\ude -v)(z_0 ) > (\chi (\de )/2 -1)\de ^{\al }.
$$
The proposition follows by comparing this inequality with
\thetag{2.6}.
\enddemo
Applying the proposition one can find $c_7$ such that
$$
z_0 \in U=\{ v<\ude -c_7 \} \subset E .
$$
By the comparison principle [K2] and \thetag{2.4}
$$
0<\int _U (dd^c \ude +\om )^n \leq \int _U (dd^c v +\om )^n\leq
\int _E (dd^c v +\om )^n =\int _E g\om ^n =0.
$$
This contradiction shows that the choice of small enough $\de$
with \linebreak $\chi(\de ) >\max (9, \chi (N\de )$ is impossible.
Thus the proof of the claim and that of the  theorem is completed.

\bigskip
\bf\centerline{References}\bigskip \rm

\noindent\item {[BT1] }{E. Bedford and B.A. Taylor, \it
 The Dirichlet problem for the complex Monge-Amp\`ere operator, \rm Invent. Math. \bf 37
\rm (1976), 1-44.}

\noindent\item {[BT2] }{E. Bedford and B.A. Taylor, \it  A new
capacity for plurisubharmonic
 functions, \rm
Acta Math. \bf 149 \rm (1982), 1-40.}

\noindent\item {[D] }{J.-P. Demailly, \it Regularization of closed
positive currents and intersection theory, \rm J. Alg. Geom. \bf 1
\rm  (1992), 361-409.}

\noindent\item {[GKZ] }{V. Guedj, S. Ko\l odziej, A. Zeriahi, \it
H\"older continuous solutions to the complex  \MA\ equations, \rm
 math.CV/0607314. }

\noindent\item {[EGZ] }{ P. Eysssidieux, V. Guedj, A. Zeriahi, \it
Singular K\"ahler-Einstein metrics, \rm \linebreak
math.AG/0603431. }

\noindent\item{[K1] }{S. Ko\l odziej, \it The complex \MA
equation, \rm Acta Math. \bf 180 \rm (1998), 69-117.}

\noindent\item {[K2] }{S. Ko\l odziej, \it   Stability of
solutions to the complex Monge-Amp\`ere on compact K\"ahler
manifolds,  \rm Indiana U. Math. J. \bf 52 \rm (2003), 667-686.}

\noindent\item{[K3] }{S. Ko\l odziej, \it The complex \MA equation
and pluripotential theory, \rm Memoirs of AMS, \bf 840 \rm (2005),
pp. 62.}

\noindent\item {[ST] }{J. Song and G. Tian, \it The K\"ahler-Ricci
flow on surfaces of positive Kodaira dimension, \rm
math.DG/0602150 .}

\noindent\item {[TZ] }{G. Tian and Z. Zhang, \it On the
K\"ahler-Ricci flow on projective manifolds of general type, \rm
Chinese Ann. Math. B \bf 27 (2) \rm (2006), \rm 179-192.}

\noindent\item {[Y] }{S.-T. Yau, \it On the Ricci curvature of a
compact K\"ahler manifold and  the complex Monge-Amp\`ere
equation, \rm Comm. Pure and Appl. Math. \bf 31 \rm (1978),
339-411.}

\noindent\item {[Z] }{Z. Zhang, \it On Degenerated Monge-Ampere
Equations over Closed K\"ahler Manifolds, \rm IMRN, vol. 2006,
1-18.}

\vskip 0.5cm

postal address: \+Jagiellonian University, Institute of
Mathematics\cr \+Reymonta 4, 30-059 Krak\'ow, Poland\cr \+e-mail:
Slawomir.Kolodziej\@ im.uj.edu.pl\cr

\bye